\documentclass{amsart}
\usepackage{graphicx} % Required for inserting images
\usepackage{color}

\newcommand{\pp}{\partial}

\newcommand{\dist}{\operatorname{dist}}
\newcommand{\ep}{\varepsilon}

\numberwithin{equation}{section}

\newtheorem{theorem}{Theorem}

\newtheorem{proposition}{Proposition}
\newtheorem{lemma}{Lemma}
\newtheorem{corollary}{Corollary}

\theoremstyle{remark}
\newtheorem{remark}{Remark}

\numberwithin{theorem}{section}
\numberwithin{proposition}{section}
\numberwithin{lemma}{section}
\numberwithin{corollary}{section}
\numberwithin{numb}{section}

\newcommand{\diam}{\operatorname{diam}}
\newcommand{\fom}{f_\Omega}
\newcommand{\dom}{d_{\Omega}}

\newcommand{\R}{\mathbb{R}}

\newcommand{\calH}{\mathcal{H}}

\newcommand{\calL}{\mathcal{L}}	
	
\newcommand{\Cconz}{{C^{con}_0}(\bar \Omega)}
\newcommand{\oO}{\omega_\Omega}

\title{Alexandrov's estimate revisited}

\date{\today}

\author{Charles Griffin}
\address{Department of Mathematics, University of Toronto, 
			Toronto, Ontario M5S 2E4, Canada.}
\email{charlie.griffin@mail.utoronto.ca}

\author{Kennedy Obinna Idu}
\address{Department of Mathematics, University of Toronto, 
			Toronto, Ontario M5S 2E4, Canada.}
\email{o.idu@utoronto.ca}

\author{Robert L. Jerrard}
\address{Department of Mathematics, University of Toronto, 
			Toronto, Ontario M5S 2E4, Canada.}
\email{rjerrard@math.toronto.edu}


\begin{document}
\maketitle

\begin{abstract}
Alexandrov's estimate states that if $\Omega$ is a bounded open convex domain in $\R^n$ and $u:\bar \Omega\to \R$ is a convex solution of the Monge-Amp\`ere equation $\det D^2 u = f$ that vanishes on $\pp \Omega$, then 
\[
|u(x) - u(y)| \le \omega(|x-y|)(\int_\Omega f)^{1/n} 
\qquad \mbox{for }\omega(\delta) = C_n\,\mbox{diam}(\Omega)^{\frac{n-1}n} \delta^{1/n}.
\]
We establish a variety of improvements of this, depending on the geometry of $\pp \Omega$.
For example, we show that if the curvature is bounded away from $0$, then the estimate remains valid if $\omega(\delta)$ is replaced by $C_\Omega \delta^{\frac 12 + \frac 1{2n}}$. We determine the sharp constant $C_\Omega$ when $n=2$, and when $n\ge 3$ and $\pp \Omega$ is $C^2$, we determine the sharp asymptotics of the optimal modulus of continuity $\oO(\delta)$ as $\delta\to 0$.  For arbitrary convex domains, we characterize the scaling of the optimal modulus $\oO$. 
Our results imply in particular that unless $\pp \Omega$ has a flat spot, $\oO(\delta) =  o(\delta^{1/n})$ as $\delta\to 0$, and 
under very mild nondegeneracy conditions, they yield the improved H\"older estimate, $\oO(\delta) \le C \delta^\alpha$ for some $\alpha>1/n$.

\end{abstract}

\section{Introduction}

Alexandrov's estimate states that if $\Omega$ is a bounded open convex domain in $\R^n$, and $u:\bar \Omega\to \R$ is a convex function such that $u=0$ on $\pp\Omega$, then there exists a constant $C_n$ such that
\begin{equation}\label{alex1}
[u]_{1/n} \le C_n \,\mbox{diam}(\Omega)^{\frac{n-1}n}  |\pp u(\Omega)|^{\frac 1n}.
\end{equation}
Here
\begin{equation}\label{basic}
    [u]_\alpha := \sup_{x,y\in \Omega, x\ne y
} \frac{|u(x)-u(y)|}{|x-y|^\alpha},
\end{equation}
and $\pp u$ denotes the subgradient of $u$, whose definition is recalled in \eqref{subgradient2}. For now we just mention that if $u$ is $C^2$, then
$|\pp u(\Omega)| = \|\! \det D^2 u\|_{L^1(\Omega)}$. %The scaling of \eqref{alex1} and its improvements below is consistent in that $|\pp (\lambda u)(\Omega)|^{\frac 1n} = \lambda |\pp u(\Omega)|^{\frac 1n}$ for every $\lambda>0$.

Estimate \eqref{alex1} plays an important role in the regularity theory of the Monge-Amp\`ere equation, see for example \cite{Gutierrez, Figalli-book}, and it is a key ingredient in some basic linear elliptic PDE estimates, see for example \cite{GT}, Chapter 9.

In this paper we establish some improvements of \eqref{alex1}. Before stating them we introduce some notation. We will write
\[
\Cconz := \{ u\in C(\bar \Omega) : \ u\mbox{ is convex}, u=0\mbox{ on }\pp \Omega\}
\]
and
\begin{equation}\label{oO.def}
\oO(\delta):= \sup\left\{ \frac{ |u(x) - u(y)|}{|\pp u(\Omega)|^{1/n}} :  u\in \Cconz, \, u\mbox{ nonzero, }\,  %\ x,y\in \bar \Omega, 
\ |x-y|\le \delta \right\}.
\end{equation}
The definition immediately implies that  for every $u\in\Cconz$, 
\begin{equation}\label{oO.est}
[ u]_{\oO} := 
\sup_{x,y\in \Omega, x\ne y }\frac{|u(x) - u(y)|}{\oO(|x-y|) } \le |\pp u(\Omega)| ^{1/n}
%\qquad\qquad \mbox{ for all }x,y\in\bar\Omega
\end{equation}
and that this is sharp in that it fails for some $u\in \Cconz$ if $\oO$ is replaced by any smaller function.
With this notation, Alexandrov's estimate \eqref{alex1} amounts to the assertion that $\oO(\delta)\le C(\Omega)\delta^{1/n}$ for all $\delta>0$. 

In this paper, we give a precise description of $\oO$, depending on the geometry of $\pp \Omega$. This allows us to show that for any bounded, convex domain $\Omega$ whose boundary satisfies a very weak nondegeneracy condition (see \eqref{t2.equiv3}), there exists some $\alpha>1/n$ such that $\oO(\delta) \le C(\Omega,\alpha)\delta^\alpha$ for all $\delta>0$, or in other words,  that
\begin{equation}\label{alex.A}
[u]_\alpha \le C(\alpha,\Omega)|\pp u(\Omega)|^{\frac 1n} \qquad\mbox{ for all }u\in \Cconz.
\end{equation}
Beyond that, we aim to characterize the range of $\alpha$ for which an estimate like the above holds, and to estimate the optimal constant $C(\alpha,\Omega)$ in \eqref{alex.A}, in terms of the geometry of $\pp\Omega$.

Our first result addresses domains for which the Gaussian curvature $\kappa$ of the boundary satisfies
\begin{equation}\label{kmin}
\inf_{\pp \Omega}\kappa = \kappa_0 >0.
\end{equation}
Except where stated otherwise, we do not impose any smoothness conditions beyond those that follow from convexity, which imply that $\pp \Omega$ is twice differentiable, and hence the Gaussian curvature is defined, $\calH^{n-1}$ a.e.. The left-hand side of  \eqref{kmin} should be understood to mean the infimum over all points at which $\kappa$ is defined.

\begin{theorem}\label{T1}
Assume that $\Omega\subset \R^n$ is convex and bounded, and that \eqref{kmin} holds.
Let 
$\alpha_* := \frac 12 + \frac 1{2n}$.
Then
%\[
% \sup_{\delta>0} \frac{\oO(\delta)}{\delta^\alpha} =:  C(\alpha,\Omega) \ \  \mbox{ is finite} 
%\qquad\mbox{ if and only if }\qquad 0<\alpha\le \alpha_* = \frac 12+ \frac 1{2n}.
%\]
%Moreover, 
\begin{align}\label{kpn2}
\mbox{  if $n=2$ then 
} \sup_{\delta>0} \frac{\oO(\delta)}{\delta^{\alpha_*}}  &=  \left( \frac{2^{3/2} }{\pi \sqrt{\kappa_0}  }\right)^{1/2} \ ,\\
\label{kpnge3}
\mbox{  If $n\ge 3$ and $\pp\Omega$ is $C^2$, then 
}
\lim_{\delta \searrow 0} \frac{\oO(\delta)}{\delta^{\alpha_*} } &=  
 \left(\frac{2^{(n+1)/2}}{|B^n_1| \sqrt{\kappa_0}}
\right)^{1/n}
\end{align}
where $|B^n_1|$ denotes the volume of the unit ball in $\R^n$.

\end{theorem}

\begin{remark}\label{rem:1}
The theorem implies that for $n=2$, 
\[
[u]_{3/4} \le \left( \frac{2^{3/2} }{\pi \sqrt{\kappa_0}  }\right)^{1/2} |\pp u(\Omega)|^{1/2} \qquad\mbox{ for all }u\in \Cconz,
\]
and that the estimate is sharp in the sense that it does not hold for any larger H\"older exponent or any smaller constant. Similarly, for $n\ge 3$, since $\oO(\delta)$ is continuous (this follows from \eqref{alex1} and the subadditivity of $\oO$, which is easily deduced from the definition) and constant for $\delta > \operatorname{diam}(\Omega)$, the theorem implies that
$\sup_{\delta>0} \delta^{-\alpha}\oO(\delta) <\infty$, and hence that \eqref{alex.A} holds, if and only if $\alpha \le \alpha_*$. 

Note also that conclusion \eqref{kpnge3} may be described as an
asymptotically sharp bound 
%of $|\pp u(\Omega)|^{1/n}( \left(\frac{2^{(n+1)/2}}{|B^n_1| \sqrt{\kappa_0}} \right)^{1/n}+o(1))$  
for the H\"older-$\alpha_*$ constant of $u\in \Cconz$ on scales $\le \delta$, as $\delta\to 0$. 
It is natural to ask
\[
\mbox{for $n\ge 3$,  is it true that }\sup_{\delta>0}\frac{\oO(\delta) }{ \delta^{\alpha_*}} =  \left(\frac{2^{(n+1)/2}}{|B^n_1| \sqrt{\kappa_0}}
\right)^{1/n}
? 
\]
This again would yield the sharp constant in \eqref{alex.A} for the critical space.
We are tempted to conjecture that the answer is ``yes", but we do not have any evidence to support this. We believe that the requirement that $\Omega$ is $C^2$ is unnecessary, and that convexity and \eqref{kmin} should suffice for \eqref{kpnge3}.
\end{remark}

Our other main result is less precise but completely general,  in particular applying to domains for which the boundary curvature may vanish. As we will see, it implies that we can improve \eqref{alex1} to stronger H\"older norms as long as the domain satisfies a very weak nondegeneracy condition. It requires more notation.
If $\Omega\subset \R^n$ is a convex set, %the support function $\sigma_\Omega$ is defined by
%\[
%\sigma_\Omega(y) = \sup_{x\in \Omega}x\cdot y.
%\]
we write $\Omega^\circ$ to denote the {\em polar} of $\Omega$, defined
by
\[
\begin{aligned}
\Omega^\circ &:= \{ y\in \R^n : x\cdot y \le 1 \mbox{ for all }x\in \Omega\} %\\
%& = \{ y\in \R^n : \sigma_\Omega(y)\le 1\}.
\end{aligned}
\]
For $a\in \Omega$ and $\nu\in S^{n-1}$, we write
\begin{align*}
S(a,\nu)& := \{ x\in \R^n : x\cdot \nu = 0, \ \ a+x\in \Omega\} \\%= (\Omega-a)_{\nu^\perp},\\
S^\circ(a,\nu) &:=  \{ y\in \R^n : y\cdot \nu = 0, \ \ x\cdot y\le 1 \ \ \mbox{ for all }x\in S(a,\nu)\}\\
&=\mbox{
 polar of $S(a,\nu)$ {\em within the hyperplane} $\nu^\perp = \{ x\in \R^n : x\cdot\nu=0\}$.}
%(\Omega-a)_{\nu^\perp}^\circ \subset \nu^\perp.
\end{align*}
If $P\subset \R^n$ is a $k$-dimensional subspace and  $A\subset P$
is a subset with relatively open interior, we will write
\[
|A| := \calH^k(A) = \mbox{$k$-dimensional Hausdorff measure of $A$} .
\]
For example,
\[
|\Omega^\circ| = \calL^n(\Omega^\circ), \qquad \qquad |S^\circ(x,\nu)| := \calH^{n-1}(S^\circ(x,\nu)), \qquad\mbox{etc.}
\]
as long as $S^\circ(x,\nu)$ has open interior within $\nu^\perp$,  which will always be the case for us.

%We also write $\pi_P:\R^n\to P$ to denote orthogonal projection onto $P$.

For $a\in \Omega$, we will write
\begin{equation}\label{dom.def}
\dom(a) := \dist(a,\pp\Omega) = \min_{b\in \pp\Omega}|a-b|
\end{equation}
and
\begin{equation}\label{Nofa.def}
N(a) = \left\{ \nu\in S^{n-1} : \exists y\in \pp \Omega \mbox{ such that }|a-y|=\dom(a)\mbox{ and }\nu = \frac{y-a}{|y-a|} \right\}    
\end{equation}
for the set of outer unit normals at boundary points closest to $a$.

We now state our second main result.

\begin{theorem}\label{T2}
Assume that $\Omega$ is a bounded, convex and open subset of $\R^n$.
Then for every positive $\delta\le \max_{a\in \Omega}\dom(a)$, 
\begin{equation}\label{oO.estB}
%\oO(\delta) \le \sup\left\{ \left( \frac{\delta}{|S^\circ(a,\nu)|}\right)^{1/n} : a\in \Omega, \dom(a) = \delta, \nu\in N(a) \right\}
 \sup_{\dom(a)=\delta} \ \sup_{\nu\in N(a)} \left( \frac{ \delta}{ 2|S^\circ(a,
\nu)|}\right)^{1/n}
\le
\oO(\delta)
 \le \sup_{\dom(a)=\delta} \ \inf_{\nu\in N(a)} \left( \frac{ n\delta}{|S^\circ(a,
\nu)|}\right)^{1/n} .
\end{equation}
\end{theorem}

In principle,  given a domain $\Omega$ with vanishing curvature, estimate \eqref{oO.estB} allows us to determine the exact scaling of $\oO(\delta)$ as $\delta\searrow 0$, and hence 
the exact range of exponents $\alpha>1/n$ for which estimate \eqref{alex.A} holds.
We illustrate this in Section \ref{sec:examples} below with several examples. For now we note

\begin{corollary}
For $\Omega$ as above and $\alpha>1/n$,  the H\"older-$\alpha$ estimate \eqref{alex.A} holds if and only if
%\item
%there exist positive constants $\beta$ and $c_\beta$ such that
%\begin{equation}\label{t2.equiv2}
%|S^\circ(a,\nu)| \ge c_\beta \dom(a)^{-\beta} \qquad\mbox{ for all }a\in \Omega.
%\end{equation}
\begin{equation}\label{t2.equiv3}
\liminf_{\delta\to 0} \inf_{\dom(a)=\delta}  \sup_{\nu\in N(a)}\delta^{\beta} |S^\circ(a,\nu)| >0\, 
\end{equation}
for $\beta = n\alpha-1$.
\end{corollary}
We omit the proof, as this follows directly from Theorem \ref{T2}.

\begin{remark}
It is known that if $S\subset\R^n$ is any bounded convex set with nonempty interior, then 
\begin{equation}\label{mahler}
|S| \ |S^\circ| \ge c_n, 
\end{equation}
see \cite{Kuperberg} for a proof with a good estimate of $c_n$ (whose sharp value is the focus of the Mahler conjecture). Thus \eqref{oO.estB} implies that there exists $C = C_n$ such that
\begin{equation}\label{oO.estbprime}
\oO(\delta)  \le  C \delta^{1/n} \sup_{\dom(a)=\delta} \ \inf_{\nu\in N(a)}  |S(a,\nu) |^{1/n} .
\end{equation}
%If in addition $S$ is centrally symmetric (that is, symmetric with respect to reflections through the origin), then
%\[
%    |S|\, |S^\circ| \le C_n = |B^n_1|^2.
%\]
%The latter inequality is sharp, with equality holding whenever $S$ is a centrally symmetric ellipsoid.
\end{remark}

\begin{remark}
It is not hard to check that if $x\in \pp\Omega$ is a point at which $\pp \Omega$ is twice
differentiable, with Gaussian curvature $\kappa$, and if $\nu$ is the outer unit normal at $x$, then
\begin{equation}\label{section.kappa}
|S^\circ(x-\delta\nu, \nu)| = \frac{\sqrt \kappa |B^{n-1}_1|}{(2\delta)^{(n-1)/2}} (1+o(1)) \qquad\mbox{ as }\delta\to 0.
\end{equation}
We present the short proof in Lemma \ref{lem:ellipse.section}.
%From this and Theorem \ref{T2}, a straightforward argument shows that under the assumptions of Theorem \ref{T1} for $n\ge 3$,
%\[
%\left(\frac {2^{(n-1)/2}}{2\sqrt{\kappa_0}|B^{n-1}_1|}\right)^{\frac 1n} \le %\liminf_{\delta\to 0}\frac{\oO(\delta)}{\delta^{\alpha_*}} \le 
%\limsup_{\delta\to 0}\frac{\oO(\delta)}{\delta^{\alpha_*}}
%\le \left(\frac {n 2^{(n-1)/2}}{\sqrt{\kappa_0}|B^{n-1}_1|}\right)^{\frac 1n}.
%\]
%This gives a quick proof of an estimate like that of Theorem \ref{T1}, though somewhat less precise.
%\end{remark}
%\begin{remark}
%    Note that \eqref{section.kappa} 
    This provides a quantitative link between the curvature at $x\in \pp\Omega$ and the rate of blowup of 
    $|S^\circ( x-\delta\nu, \nu)|$ as $\delta\searrow 0$. In view of this, it is natural to interpret  \eqref{t2.equiv3} as a {\em degenerate} positive curvature condition, growing more degenerate (and yielding a weaker H\"older exponent) as $\beta$ decreases.
\end{remark}

To conclude this introduction, we note that several recent works have established sharp estimates of H\"older seminorms of solutions of Monge-Amp\`ere equations of the form
\begin{equation}\label{MA}
\det D^2 u = F(x, u, Du)
\end{equation}
for particular geometrically meaningful functions $F(x, u, Du)$, see for example \cite{JianLi, Le22, Le23,jian2023sharp, chen2023anisotropic}. Some of these papers allow for domains in which the boundary curvature can vanish, and they determine H\"older exponents that reflect the boundary behaviour in a way that has some similarities to what we find in Theorem \ref{T2}; see Corollary \ref{cor:1}. The proofs in these references rely on careful constructions of sub- and supersolutions or even in rare cases explicit solutions. These play no role in our arguments.

\section{preliminaries, and the proof of Theorem \ref{T2}}

Like all the results in this paper, those in this section are elementary, and many if not all (apart from the proof of Theorem \ref{T2}, which however is an immediate corollary of other results) are presumably known to experts. For the convenience of the reader, we nonetheless provide complete proofs, mostly self-contained.

First, we recall some standard definitions. For $u\in \Cconz$ and $x\in \Omega$,
\begin{equation}\label{subgradient1}
\pp u(x) := \{ p\in \R^n : u(x) + p \cdot (y-x) \le u(y)\ \mbox{ for all }y\in \Omega \}
\end{equation}
and for $A\subset\Omega$,
\begin{equation}\label{subgradient2}
\pp u(A) := \cup_{x\in A} \pp u(x). % \qquad\mbox { for $A\subset \Omega$}.
\end{equation}
As mentioned above, if $u$ is $C^2$ and strictly convex, then by the change of variables $p = Du(x)$,
\[
\|\! \det D^2 u\|_{L^1(\Omega)} = \int_\Omega \det D^2 u  = \int_{p\in Du(\Omega) } \,dp = 
|\pp u(\Omega)|.
\] 
(This remains true under somewhat weaker assumptions.)
Given $a\in \Omega$ we will write $u_{a}:\bar\Omega\to \R$ to denote the function defined by
\[
u_{a}((1-\theta) y + \theta a) = -\theta \qquad\mbox{ for every }y\in \pp \Omega\mbox{ and }\theta\in [0,1]. 
\]
The definition states that
\[
u_a=0\mbox{ on }\pp\Omega, \qquad u_a(a)=-1,
\]
and $u_a$ is linear on the line segment from any point on $\pp\Omega$ to $a$.
When we wish to explicitly indicate the dependence of $u_a$ on $\Omega$, we will write $u_{\Omega,a}$.
It is well-known and straightforward to check that $u_a$ is convex.

Next, we define $\fom:\Omega\to \R$ by
\[
\fom(a) := |\pp u_{a}(\Omega)| = \calL^n(\pp u_{a}(\Omega))\qquad\quad\mbox{ where }u_a = u_{\Omega,a}.
\]
The following result implies that to understand the modulus of continuity for functions $u\in \Cconz$ with $|\pp u(\Omega)|$ finite, it suffices to study the asymptotics of $\fom(a)$ as $a\to \pp\Omega$.

\begin{proposition}\label{prop:fom}
Let $\Omega$ be a bounded, convex, open subset of $\R^n$. Then the modulus  
$\oO$ defined in \eqref{oO.def} satisfies
\begin{equation}\label{oO.fom}
\oO(\delta) = \sup \{ \fom(a)^{-1/n} : \dom(a)\le \delta\}. 
\end{equation}
\end{proposition}

\begin{proof}
{\em Step 1}. We first claim that for $u\in \Cconz$ and any $a,b\in \Omega$, there exists $\bar a\in \Omega$
such that 
\begin{equation}\label{moc.claim1}
\dom(\bar a)\le |a-b|, \qquad |u(\bar a)|\ge |u(a)-u(b)|.
\end{equation}
We recall the proof, which is standard.
Consider $a,b\in \Omega$ such that $u(a) \le u(b)\le 0$.
Let $\bar b$ be the point in $\partial \Omega$ on the ray that starts at $a$ and passes through $b$. 
Then there exists some $\theta\in (0,1)$ such that
$
b = (1-\theta)a + \theta \bar b$.
We next define $\bar a =  \theta a + (1-\theta) \bar b$.
These definitions imply that
\[
\bar a - \bar b = \theta(a - \bar b) = a - b.
\]
Thus $\dom(\bar a)\le |\bar a -\bar b|  = |a-b|$.
Moreover, by convexity,
\begin{align*}
u( b) - u(a)  &= u((1-\theta)a + \theta \bar b) - u(a) \ \le \  \theta (u(\bar b) - u(a))\\
&= u(\bar b)  - [ (1-\theta)u(\bar b) + \theta u(a)]\\
&\le u(\bar b) - u((1-\theta) \bar b + \theta a) \ = \ u(\bar b) - u(\bar a) = |u(\bar a)|.
\end{align*}
Since $\dom(\bar a)\le |a - b|$, this proves \eqref{moc.claim1}. 
%It  follows that
%\begin{equation}\label{moc.c2}
%\oO(\delta)\le \sup\{ |u(a)|  \, : \,  u\in \Cconz, u\mbox{ convex},\ a\in \Omega, \dom(a) \le \delta\}.
%\end{equation}

{\em Step 2}. 
Given $u\in \Cconz$, $\delta>0$ and $x,y\in \Omega$ such that $|x-y|\le \delta$, fix $a\in \Omega$ such that $\dom(a)\le \delta$ and $|u(x)-u(y)|\le |u(a)|$, and define $w(x) = u(a) u_a(x)$.
Then $u\le w\le 0$ in $\Omega$ and $u=w=0$ on $\pp\Omega$, so standard arguments (see e.g., \cite[Lemma~1.4.1]{Gutierrez}) imply that
\begin{equation}\label{alex.mod}
|\pp u(\Omega)| \ge |\pp w(\Omega)| = |u(a)|^n \fom(a),
\end{equation}
with equality if and only if $u = w$. The definition \eqref{oO.def} of $\oO$ then implies that
\begin{equation}\label{mod.01}
|u(x)-u(y)|\le |u(a)| \overset{\eqref{alex.mod}}\le \fom(a)^{-1/n} |\pp u(\Omega)|^{1/n} %\le \oO(\delta) |\pp u(\Omega)|^{1/n}.
\end{equation}
Thus for nonzero $u\in \Cconz$,
\[
\mbox{ if }|x-y|\le \delta, \mbox{ then }\ \quad \ \frac{|u(x)-u(y)|}{ |\pp u(\Omega)|^{1/n} } \le \sup_{\dom(a)\le \delta} \fom(a)^{-1/n}.
\]
It follows from this and the definition of $\oO$ that 
\[
\oO(\delta)\le\sup_{\dom(a)\le \delta} \fom(a)^{-1/n}.
\]
On the other hand, given any $a\in \Omega$ such that $\dom(a)\le \delta$, consider $u = u_a$, and fix $b\in \pp\Omega$ such that $\dom(a) = |a-b|$. Then
\[
|u_a(a)-u_a(b)| = |u_a(a)| = 1 = \frac {|\pp u_a(\Omega)|^{1/n}}{\fom(a)^{1/n}},
\]
and thus
\[
\oO(\delta) \ge \sup_{|x-y|\le \delta} \frac{|u_a(x)-u_a(y)|}{|\pp u_a(\Omega)|^{1/n}}
\ge \fom(a)^{-1/n} \qquad\mbox{ whenever }\ \ \dom(a)\le \delta. 
\]

\end{proof}

Motivated by Proposition \ref{prop:fom}, we record some properties of $\fom$ and related notions.

\begin{lemma}\label{anywhere}
$\pp u_a(a) = \pp u_a(\Omega)$.
\end{lemma}

\begin{proof}
It is clear that $\pp u_a(a) \subset \pp u_a(\Omega)$. To prove the other inclusion, assume that 
$p\in \pp u(x_0)$ for some $x_0\in \Omega$. We must show that $p\in \pp u(a)$. We may assume that
$x_0\ne a$, so we can write $x_0 = \theta a + (1-\theta)y$ for some $y\in \pp \Omega$ and $\theta\in (0,1)$.

For $x\in \Omega$ and $p\in \R^n$ we will write $\ell_{x,p}(z) := u_a(x) + p\cdot (z-x)$, so that $p\in \pp u_a(x)$ if and only if $\ell_{x,p}\le u_a$ in $\Omega$. 
%The definition of subgradient implies that $ u_a(z)\ge \ell_{x_0,p}(z)$ for all $z\in \Omega$. Evaluating this at $z = s a + (1-s)y$  for $s\in (0,1]$ and using the definition of $u_a$, we find that
%\[
%-s \ge -\theta + (s-\theta) p\cdot (a-y)\qquad\mbox{ for all }s\in (0,1].
%\]
%Since $s-\theta$ can be either positive or negative as $s$ varies, this implies that $p\cdot (a-y) = 1$. It follows that $\ell_{x_0, p} = \ell_{a,p}$, since
%\[
%\ell_{x_0, p}(z) - \ell_{a,p}(z) =  1-\theta + p\cdot (a - x_0) = 1-\theta + p\cdot [(1-\theta)(a-y)] = 0 \quad\mbox{ for all }z.
%\]
%Thus $\ell_{a,p} = \ell_{x_0,p} \le u_a$ in $\Omega$, so $p\in \pp u_a(a)$.
Since $u_a$ and $\ell_{x_0, p}$ are both linear when restricted to the segment
$\{ sa+(1-s)y : s\in (0,1]\}$, and because $x_0$ belongs to the interior of this segment and $u_a\ge \ell_{x_0,p}$ on this segment, we see that $u_a = \ell_{x_0,p}$ on this segment, and in particular at $x = a$. Thus $\ell_{x_0,p}$ is a supporting hyperplane at $a$; in fact $\ell_{x_0,p} = \ell_{a,p}$. It follows that $p\in \pp u_a(a)$.
\end{proof}

\begin{lemma}\label{lem:compare}
If $a \in \Omega\subset \Omega'$ then $f_{\Omega}(a) \ge f_{\Omega'}(a)$.
\end{lemma}

\begin{proof}
If $p\in \pp u_{\Omega',a}(\Omega')$ then $p\in \pp u_{\Omega', a}(a)$, which implies that
$\ell_{a,p}\le u_{\Omega',a}$ in $\Omega'$. But it is easy to check that 
$u_{\Omega',a} \le u_{\Omega,a}$ in $\Omega$, and it follows that $\ell_{a,p}\le u_{\Omega,a}$
in $\Omega$, which implies that $p\in \pp u_{\Omega, a}(a)$.

Thus
$\pp u_{\Omega',a}(\Omega') \subset \pp u_{\Omega,a}(\Omega)$, from which we deduce that
$f_{\Omega'}(a)\le f_\Omega(a).$
\end{proof}

\begin{lemma}
Assume that $\Omega\subset\R^n$ is bounded, convex and open, with $a\in \Omega$.

Then
\[
\fom(a)= |(\Omega-a)^\circ|, \qquad\mbox{ where }\Omega-a = \{ x-a : x\in \Omega\}.
\]
\end{lemma}

An equivalent statement appears as an exercise (problem 3.3) in the recent text \cite{NamLe-book}.

\begin{proof}
We first prove the lemma for $a=0$.
We know from Lemma \ref{anywhere} that $\pp u_0(\Omega) = \pp u_0(0)$. Then
\begin{align*}
p\in \pp u_0(0) 
&\iff  u_0(x) \ge u_0(0) + p\cdot x\mbox{ for all }x\in \Omega\\
&\iff  u_0(x) \ge -1 + p\cdot x\mbox{ for all }x\in\bar \Omega
\end{align*}
Since $u_0(x)\le 0$ in $\Omega$, it follows that 
\[
p\in \pp u_0(0) \ \ \Longrightarrow \ \ -1+p\cdot x\le 0\mbox{ for all }x\in \Omega \ \  \Longrightarrow \ \  p\in \Omega^\circ.
\]
On the other hand, if $p\in \Omega^\circ$ then $\ell_{0,p}(x):= -1+x\cdot p$
is an affine function such that  $\ell_{0,p} \le 0 = u_0$ on $\pp \Omega$ and $\ell_{0,p}(0)=-1 = u_0(0)$. It follows from this and the definition of $u_0$ that $\ell_{0,p} \le u_0$ in $\Omega$, and hence that $p\in \pp u_0(0)$.

It follows that $\pp u_0(\Omega) = \Omega^\circ$, and hence that $\fom(0) = |\Omega^\circ|$.

\smallskip

For general $a\in \Omega$, the definitions imply that for every $x\in \Omega$, 
\[
u_{\Omega, a}(x) = u_{\Omega-a,0}(x-a), \qquad\mbox{ and thus }\pp u_{\Omega,a}(x) = \pp u_{\Omega-a,0}(x-a).
\]
Thus $\fom(a) = |\pp u_{\Omega, a}( a) | = |\pp u_{\Omega-a,0}(0)| = |(\Omega-a)^\circ| $.
\end{proof}

% \begin{lemma}
% $\sigma_{\Omega-a}(y) = \sigma_\Omega(y) - a \cdot y$.
% \end{lemma}

% \begin{proof}
% \[
% \sigma_{\Omega-a}(y) = 
% \sup \{ x\cdot y : x\in \Omega-a\} = \
% \sup \{ (z-a)\cdot y : z\in \Omega\} = \sigma_\Omega(y)-a\cdot y.
% \]
% \end{proof}

\begin{lemma}\label{lem:affine}
Let $M:\R^n\to \R^n$ be an invertible linear transformation, and let
$M\Omega := \{ Mx: x\in \Omega\}$.
Then
\[
f_{M\Omega}(Ma) = |\!\det M|^{-1} \fom(a).
\]
\end{lemma}

\begin{proof}
The definitions imply that for every $x\in \Omega$,
\[
u_{M\Omega,  Ma}(Mx) = u_{\Omega,a}(x).
\]
Thus for every $x\in \Omega$,
\begin{align*}
&p\in \pp u_{M\Omega, Ma}(Mx) \\
&\qquad \ \  \Longleftrightarrow \ \ 
u_{M\Omega, Ma}(Mz) \ge u_{M\Omega, Ma}(Mx) + p \cdot (Mz-Mx)\qquad\mbox{ for all }Mz\in M\Omega\\
&\qquad \ \  \Longleftrightarrow \ \ 
u_{\Omega, a}(z) \ge u_{\Omega, a}(x) + M^Tp \cdot (z-x) \qquad\mbox{ for all }z\in \Omega.
\\
&\qquad \ \  \Longleftrightarrow \ \ 
M^Tp\in \pp u_{\Omega,a}(x).
\end{align*}
We deduce that $\pp u_{M\Omega, Ma}(M\Omega) = \{ M^{-T}p : p\in \pp u_{\Omega,a}(\Omega)\}$.
Now the conclusion follows from basic properties of Lebesgue measure.
\end{proof}

\begin{lemma}\label{lem:slice}
Assume that $\Omega\subset \R^n$ is a bounded, open convex set containing the origin. %and that
For any subspace $P$ of $\R^n$, define
\[
\Omega_P := \Omega\cap P, \qquad\qquad \Omega_P^{\circ} := \{ y\in P : x\cdot y\le 1 \mbox{ for all }x\in \Omega_P\}.
\]
(Thus $\Omega_P^\circ$ denotes the polar of $\Omega_P$ within $P$ rather than within the ambient $\R^n$.)

Let $\pi_P:\R^n\to P$ denote orthogonal projection onto $P$.

Then
\[
\Omega_P^\circ = \pi_P(\Omega^\circ).
\]
\end{lemma}

\begin{proof}
%We have
%\begin{align*}
%y\in \Omega^\circ 
%&\ \ \Longrightarrow \ \  y\cdot x\le 1\quad\mbox{ for all }x\in \Omega\\
%&\ \ \Longrightarrow \ \  y\cdot x\le 1\quad\mbox{ for all }x\in \Omega_P \\
%&\ \ \Longrightarrow \ \  \pi_P y\cdot x\le 1\quad\mbox{ for all }x\in \Omega_P  \ \ \Longrightarrow \ \  \pi_Py\in \Omega_P^\circ,
%\end{align*}
%proving that $\pi_P(\Omega^\circ) \subset \Omega_P^\circ$.

%To establish the converse, w
We will show that $(\pi_P(\Omega^\circ))^\circ = (\Omega_P^\circ)^\circ = \bar\Omega_P$, where our convention is that if $A\subset P$ is convex, then $A^\circ $ denotes the polar within $P$, whereas if $A$ is a convex set not contained in $P$, then $A^\circ$ denotes its polar in $\R^n$.
Then
\begin{align*}
y\in (\pi_P(\Omega^\circ))^\circ
&\ \ \Longleftrightarrow \ \  y\in P\mbox{ and }y\cdot x\le 1\quad\mbox{ for all }x\in \pi_P(\Omega^\circ)\\
&\ \ \Longleftrightarrow \ \  y\in P\mbox{ and } y\cdot \pi_P x\le 1\quad\mbox{ for all }x\in \Omega^\circ\\
&\ \ \Longleftrightarrow \ \   y\in P\mbox{ and } y\cdot x\le 1\quad\mbox{ for all }x\in \Omega^\circ\\
&\ \ \Longleftrightarrow \ \ y\in P \cap (\Omega^\circ)^\circ = P\cap \bar \Omega = \bar \Omega_P,
\end{align*}
completing the proof.

%To establish the converse, fix a nonzero $y_0\in \Omega_P^\circ$, and let $x_0$ maximize $x\mapsto y_0\cdot x$ in $\bar \Omega$. Since $\Omega$ is bounded, such a maximizer exists, and clearly it must belong to $\pp \Omega$. It follows that $x_0$ satisfies the Lagrange multiplier equations, which may be written
%\[
%x_0\in \pp \Omega_P, \qquad y_0 = \lambda \nu_{\Omega_P}(x_0) \mbox{ for some }\lambda>0
%\]
%where $\nu_{\Omega_P}$ is an outer unit normal  to $\pp\Omega_P$ within the plane $P$.

%It follows from the Hahn-Banach Theorem, if you like, that there exists an outer unit normal $\nu_\Omega(x_0)$ to $\pp \Omega$ in $\R^n$ such that $\pi_P \nu_\Omega(x_0) = a \nu_{\Omega_P}(x_0)$ for some $a>0$. Let $Y_0 = \frac \lambda a \nu_\Omega(x_0)$.  It is clear that $\pi_PY_0=y_0$, and we will show that $Y_0\in \Omega^\circ$ and $\pi_P Y_0 = y_0$. Indeed, this is a consequence of the fact that $x_0$ satisfies the Lagrange multiplier equations
%\[
%x_0\in \pp \Omega, \qquad Y_0 = \frac \lambda a \nu_{\Omega}(x_0) \mbox{ for some }\lambda,a>0.
%\]
%Then by the convexity of $\Omega$, it follows that $x_0$ maximizes $x\mapsto Y_0\cdot x$ in $\bar \Omega$. But since $x_0\in \Omega_P\subset P$, 
%\[
%x_0 \cdot Y_0 = x_0 \cdot \pi_P Y_0 = x_0\cdot y_0 \le 1, \quad\mbox{ because }y_0\in \Omega_P^\circ.
%\]
%Thus $x\cdot Y_0\le 1$ for all $Y_0\in \bar\Omega$, proving the claim.
\end{proof}

\begin{lemma}\label{lem:workhorse}
Let $\Omega$ be an open convex subset of $\R^n$ with nonempty boundary. For $a\in \Omega$, let $x\in \pp\Omega$
%let $t := \dom(a)$, and let $x\in \pp\Omega$ 
be a point such that $|a-x| = \dom(a)$, and let $\nu = \frac{x-a}{|x-a|}$ . 
(Thus $\nu\in N(a)$, in the notation introduced in \eqref{Nofa.def}.)
Then
\begin{equation}\label{workhorse}
\frac 1 n \dom(a)^{-1} |S^\circ(a,\nu)| \le \fom(a)  \le 2 \dom(a)^{-1} |S^\circ(a, \nu)|.
\end{equation}
\end{lemma}

\begin{remark}
A curious consequence of \eqref{workhorse} is that for $a\in \Omega$, if there exist more than one point $b\in \pp \Omega$ such that $\dom(a) = |b-a|$, then
\[
\sup_{\nu \in N(a)} |S^\circ(a,\nu)| \le 2n \inf_{\nu \in N(a)} |S^\circ(a,\nu)| 
\]
\end{remark}

\begin{proof}
{\em Step 1}. After a translation and a rotation, we may assume that $a=0$ and that $x = (0,\ldots, 0,-\delta)$ where $\delta = \dom(a)$. Then  $-e_n$ is the outer unit normal at $x$, and hence  $\Omega \subset \{ y\in \R^n : y_n > -\delta\}$. One can then quickly check that
\begin{equation}\label{Ocirc1}
\{ -s e_n : 0 \le s \le \frac 1 \delta\} \subset \Omega^\circ . 
\end{equation}
Let $P := \R^{n-1}\times \{0\}$, so that $\Omega_P = S(a,\nu)$.
It then follows from Lemma \ref{lem:slice} that 
\begin{equation}\label{Ocirc2}
 S^\circ(a,\nu) = \pi_P (\Omega^\circ).
 \end{equation}
 Now let $T$ denote the Steiner symmetrization of $\Omega^\circ$ with respect to the hyperplane $x_n=0$. Well-known properties of Steiner symmetrization imply that $|T| = |\Omega^\circ|$, that $T$ inherits the convexity of $\Omega^\circ$, and owing to \eqref{Ocirc1}, \eqref{Ocirc2}, that
 \[
 S^\circ(a,\nu)\subset T, \qquad \{ \pm \frac 1{2\delta} e_n \} \subset T
 \]
By convexity, $T$ contains the cones in $\R^n$ with base $S^\circ(a,\nu)\subset \R^{n-1}\times \{0\}$ with vertices at $\pm \frac 1 {2\delta}e_n$. Each of these cones has measure $\frac 1{2 \delta n} |S^\circ(a,\nu)|$.
We conclude that 
\[
|\Omega^\circ| = |T| \ge \frac 1{ \delta n} |S^\circ(a,\nu)| = \frac 1 n \dom(a)^{-1}|S^\circ(a,\nu)|.
\]

{\em Step 2}. Let $P^\perp = \{0^{n-1}\}\times \R$, the orthogonal complement of $P$.
Since $a=0$ and $\dom(a) \ge \delta$, it is clear that $\Omega\cap P^\perp = \Omega_{P^\perp} \supset \{0^{n-1}\}\times (-\delta,\delta)$. It easily follows that $\Omega_{P^\perp}^\circ \subset \{0^{n-1}\}\times (-\frac 1\delta, \frac 1\delta)$. In addition, Lemma \ref{lem:workhorse} implies that
\[
\pi_{P^\perp}(\Omega^\circ) \subset \Omega_{P^\perp}^\circ.
\]
It follows from these facts and \eqref{Ocirc2} that
\[
\Omega^\circ\subset S^\circ(a,\nu) \times (-\frac 1\delta, \frac 1\delta),
\]
(writing $S^0(a,\nu)$ as a subset of $\R^{n-1}$ rather than of $\R^{n-1}\times \{0\}$).
Thus
\[
|\Omega^\circ| \le |S^\circ(a,\nu)|\times \frac 2\delta =  2\dom(a)^{-1}|S^\circ(a,\nu)|.
\]
\end{proof}

\begin{proof}[Proof of Theorem \ref{T2}]
Estimate \eqref{oO.estB} follows directly from Proposition \ref{prop:fom} and Lemma \ref{lem:workhorse}.
\end{proof}

\section{Examples}\label{sec:examples}

Our first illustration of the utility of Theorem \ref{T2} addresses a class of convex sets considered in several recent papers.

\begin{corollary}\label{cor:1}
    Let $\Omega$ be a bounded, open convex subset of $\R^n$, and assume that there exist positive constants $\eta$ and $a_1,\dots, a_k$, with $k\le n-1$, such that at any $b\in \pp \Omega$, after a translation and a rotation, 
    \begin{equation}\label{skc}
b=0\quad  \text{ and }\quad \Omega\subseteq\{x\in\R^n: x_{n}>\eta (|x_1|^{p_1}+\cdots + |x_k|^{p_k})\}.     
    \end{equation}
    Then there exists a constant $C$, depending on $\eta, n, \diam(\Omega)$, such that %$\oO(\delta) \le C \delta^{\alpha}$ 
    \[
    [u]_\alpha \le C|\pp u(\Omega)|^{1/n} \ \ \ \mbox{ for all }u\in \Cconz, \qquad\mbox{ where }\alpha = \frac 1n(1 + \sum_{j=1}^k  \frac 1{p_j}).
    \]
\end{corollary}

Note that \eqref{skc} allows $\Omega$ to be completely degenerate at $b$  in $n-k-1$ directions. 

In \cite{jian2023sharp, chen2023anisotropic}, sharp H\"older estimates on domains satisfying \eqref{skc} at every  $b\in \pp \Omega$ (for a suitable $b$-dependent choice of coordinates) are proved for solutions of certain equations of the form \eqref{MA}. Interestingly, the quantity $\sum_{j=1}^k  \frac 1{p_j}$ also appears in the H\"older exponents in these results, modified by other parameters appearing in the nonlinearity on the right-hand side of \eqref{MA}. 

\begin{proof}
Let $\delta>0$ and $a\in \Omega$ with $d_{\Omega}(a)=|a-b|=\delta$ for some $b\in \pp\Omega$. After a translation and rotation, we may assume that \eqref{skc} holds.
We necessarily have that $a=(0, \ldots, 0, \delta)$. Indeed, suppose $a_i\ne 0$ for some $1\le i \le n-1$, then from the supporting hyperplane $\{x_{n}=0\}$, we obtain
\[
d_\Omega(a)\le \dist(a,\{x_{n}=0\})=|a_{n}|<|a|=d_\Omega(a),
\]
a contradiction, verifying the claim.

Now, relabelling coordinates, we write the unit outer normal at $b$ as $\nu=-e_{n}$ and have that
\begin{align*}
    S(a,\nu) &\subseteq \left\{x\in \R^{n-1}\times \{0\}: \delta= \eta \sum_{i=1}^k |x_i|^{p_i}, \ \  |(x_{k+1}, \ldots x_{n-1})|< 
    \mbox{diam}(\Omega)\right\}.
%    &\subseteq \left\{x\in \R^{n-1}\times \{0\}: \delta\ge\eta_i|x_i|^{p_i} \text{ for all }1\le i\le k\right\}\\
%    &= \left\{x\in \R^{n-1}\times \{0\}: |x_i|\le \left(\frac{\delta}{\eta_i}\right)^{1/p_i} \text{ for all }1\le i\le k\right\}.
\end{align*}
Thus $|S(a,\nu)|$ is bounded by the volume of the set on the right, which is
\[
C\delta^{\frac 1{p_1}+\ldots + \frac 1{p_k} } 
\]
for a constant $C$ depending\footnote{Using a formula 
derived by Dirichlet and quoted on the ``generalizations'' section of the Wikipedia page for ``Volume of an n-ball'', one can check that the constant this computation yields is
\[
C(p_1,\ldots, p_k, \eta, k, n-1) = \eta^{-(\frac 1{p_1}+\ldots + \frac 1{p_k})} 
\frac{ 2^k \Gamma(1+ \frac 1 {p_1})\cdots \Gamma(1+\frac 1 {p_k})}
{\Gamma(1+ \frac 1{p_1} +\cdots +\frac 1{p_k})} |B^{n-k-1}_1| \diam(\Omega)^{n-k-1},
\]} on $\eta, p_1,\ldots, p_k, k, n-1, \mbox{diam}(\Omega)$.
Since $a$ was arbitrary,  \eqref{mahler} and
Theorem \ref{T2} (or see \eqref{oO.estbprime}) imply that
\[
\omega_{\Omega}(\delta) \le C\,\delta^{\alpha}.
\]
Hence the result on the H\"older estimate.
\end{proof}

If there is any point $b\in \pp \Omega$ such that after a translation and a rotation
\[
b=0\quad  \text{ and }\quad \Omega\supseteq\{x\in\R^n: x_{n}>\frac 1\eta (|x_1|^{p_1}+\cdots + |x_k|^{p_k}), \ \ |(x_{k+1},\ldots, x_n )|< h\}   
\]
for some positive numbers $\eta, p_1,\ldots, p_k, h$, then by a similar argument to that above one can show that $\oO(\delta) \ge c \delta^\alpha$ for all sufficiently small $\delta$ and the same $\alpha$ as above. This would use the fact that if $S$ is a {\em centrally symmetric} convex body in $\R^k$, then 
$|S| \, |S^\circ| \le |B^k_1|^2$.

The following lemma provides a way to generate a large class of examples.

\begin{lemma}
Assume that $\Omega\subset\R^2$ is a smooth convex set of the form
\begin{equation}\label{graph.h}
\Omega = \{ (x_1,x_2)\in \R^2  : |x_1|< R,  h(x_1) < x_2 < D - h(x_1)\}
\end{equation}
for some $R, D>0$, where $h:[-R,R]\to [0,\infty)$ is an  even function, smooth on $(-R,R)$, such that $h(0)=h'(0) = h''(0)=0$ and $h(R) = \frac D2$. 

Assume moreover that the boundary curvature is nondecreasing as one moves in the direction of increasing $x_1$ along $\pp \Omega$ from $(0,  0)$ toward $(R,D/2)$.

Then, writing $h^{-1}(\delta)$ to denote the unique {\em positive} solution of the equation $h(x)=\delta$ for $0<\delta\le D/2$,  there exists $\delta_0>0$ such that 
\begin{equation}\label{2d.special}
\frac 12 \sqrt{\delta h^{-1}(\delta)} \le \oO(\delta) \le  \sqrt 2\sqrt{\delta h^{-1}(\delta)}  \qquad\mbox{ for }0<\delta <\delta_0.
\end{equation}

\end{lemma}

The lemma implies that  given any modulus of the form 
$\omega(\delta) = \sqrt{\delta h^{-1}(\delta)}$ for $h$ satisfying the above hypotheses, we can construct a domain for which the sharp modulus of continuity $\oO$ in the Alexandrov estimate exactly agrees with $\omega$, up to a factor of $2\sqrt 2$.

\begin{proof}
Assume that $0<\delta<\delta_0$, to be fixed below.

{\em Step 1.} It is clear from \eqref{graph.h} and properties of $h$ that if $\delta_0$ is sufficiently small (in fact  here  $\delta_0<D/2$ is sufficient), then the origin is the unique closest boundary point to $(0,\delta)$, and hence that $N(a)$ as defined in \eqref{Nofa.def} consists of $\{ -e_2\}$. 
Then the definitions imply that
\[
S((0,\delta), -e_2) =  ( - h^{-1}(\delta) , h^{-1}(\delta)) \times \{0\}.
\]
Recall our convention that $S^\circ(a,\nu)$ denotes the polar within the subspace $\nu^\perp$.
If $a$ and $b$ are positive numbers, then $(a,b)^\circ = [\frac 1a, \frac 1b]$, so it follows that
\[
|S^\circ((0,\delta), -e_2) | = \frac 2 { h^{-1}(\delta)}.
\]
This and Theorem \ref{T2} imply the lower bound for $\oO(\delta)$ in \eqref{2d.special}.

{\em Step 2.} To complete the proof of the Lemma, again by Theorem \ref{T2}, it suffices to show that if $\dom(a)=\delta$ and $\nu\in N(a)$, then
\begin{equation}\label{Scirc}
 |S^\circ(a,\nu)| \ge \frac 1{h^{-1}(\delta)} ,
\end{equation}
if $\delta_0$ is small enough. Fix any $a\in \Omega$ such that $\dom(a)=\delta$ and any $b\in \pp\Omega$ such that $\dom(a)=|a-b|$, and let $\nu = \frac{b-a}{|b-a|}$. Noting from \eqref{graph.h} that $\Omega$ is symmetric about the $y$ axis (since $h$ is even) and about the line $x = D/2$,  we can assume that $b\in\{(x_1, x_2)\in \pp \Omega : 0\le x_1\le R, x_2 = h(x_1)\}$.

Then we define $\widetilde \Omega$ to be the set obtained by translating $b$ to the origin and rotating so that $\widetilde \Omega\subset \{(x_1, x_2): x_2>0\}$. This operation moves $a$ to the point $(0,\delta)$. Next we let $\tilde h_1$ be the function whose graph parametrizes the lower part of $\pp \widetilde \Omega$, defined by $\tilde h(x_1) := \inf \{ x_2\in \R : (x_1, x_2)\in \widetilde\Omega\}$. 
By our assumption about the monotonicity of the boundary curvature along the short arc connecting $(0,0)$ to $(R,D/2)$, we see that if $\delta_0$ is small enough, then 
\[
\mbox{curvature of $\pp  \Omega$ at $(x_1, h(x_1))$}
\le
\mbox{curvature of $\pp \widetilde \Omega$ at $(x_1, \tilde h(x_1))$}.
%= \frac{h''(x_1)}{(1+h'(x_1)^2)^{3/2}} < \frac {\tilde h''(x_1)}{(1+\tilde h'(x_1)^2)^{3/2}} \qquad\mbox{ for }0<x_1<r_0.
\]
for $0< x_1 < h^{-1}(\delta_0)$. Since $\tilde h(0) = \tilde h'(0) = h(0) = h'(0)$, and because $0=h''(0)\le \tilde h''(0)$, this implies that $\tilde h(x_1) \ge h(x_1)$ for $0<x_1< h^{-1}(\delta_0)$.

Computing $S^\circ(a,\nu)$ in the coordinate system of $\widetilde \Omega$,  we find that 
$S(a,\nu) = (-\alpha, \beta)\times \{0\}$, where $-\alpha,\beta$ are the negative and positive solutions, respectively, of the equation $\tilde h(x)=\delta$.  and thus 
\[
|S^\circ(a, \nu)| = \frac 1\beta+ \frac 1 \alpha \ge \frac 1\beta.
\]
But the fact that $\tilde h \ge h$ for  $0<x_1< h^{-1}(\delta_0)$ implies that $\beta\le h^{-1}(\delta)$,
proving \eqref{Scirc}.
\end{proof}

Based on the above lemma, it is straightforward to construct examples of domains $\Omega\subset \R^2$ such that
$\oO(\delta) \sim \delta^{1/p}$ for given $p>2$. Another example is obtained by by taking $h(x)$ in \eqref{graph.h} such that
\[
h(0) = 0, \qquad h(x) = e^{-1/|x|}\quad\mbox{ for }0<|x|< a
\]
and extended (after choosing $a$ small enough) so that the graph of $h$ has increasing curvature until the point where its tangent becomes vertical.
Then the lemma implies that for the resulting domain $\Omega$,
\[
\frac 12 \left( \frac{\delta}{|\log \delta|}  \right)^{1/2} \le
\oO(\delta) 
\le
\sqrt 2 \left( \frac{\delta}{|\log \delta|}  \right)^{1/2}.
\]
In this spirit, it would be straightforward to construct sets with $\oO$ for example having logarithmic or other corrections to H\"older moduli $\delta^\alpha$ for some $\frac 1n <\alpha< \frac 12 + \frac 1{2n}$.

\medskip
The next lemma shows that, loosely speaking, the scaling in the classical Alexandrov estimate \eqref{alex1} is almost never optimal:

\begin{lemma}\label{lem:flatspot}
    Let $\Omega\subset \R^n$ be a bounded, convex and open domain. Then 
    \[
    \exists A, \delta_0 >0\mbox{ such that }\oO(\delta) \ge A \delta^{1/n}\mbox{ for }0<\delta<\delta_0 \qquad\Longleftrightarrow \qquad \pp \Omega \mbox{ has a flat spot}.
    \]
    In fact, if $\oO(\delta) \ge A\delta^{1/n}$ for $\delta\in (0,\delta_0)$ then there exists a supporting hyperplane $P$ such that 
    \[
    P\cap \pp\Omega  \mbox{ contains an $n-1$-dimensional ball of radius }\frac{A^n |B^{n-2}_1|}{2^{n-1} n R^{n-2}}.
    \]
\end{lemma}

The estimate of the radius of the ball is not sharp.

\begin{proof} 
We first claim that for $R, c>0$ and $S\subset \R^k$,
\begin{equation}\label{fs1} 
\mbox{ if }S\subset B_R \mbox{ and }|S^\circ|<c, \qquad\mbox{ then } B_r\subset S \mbox{ for }r= \frac {|B^{k-1}_1|}{2c(2R)^{k-1}}
\end{equation}
Indeed, for $r<R$, suppose $S\subset B_R$ does not contain $B_r$. By a rotation we may assume that there is a point of the form
$b = (0,\ldots, 0, r_1)$ with $0<r_1<r$ such that $d_S(0) = |0 - b| = r_1$. Then the plane $\{ x : x_k = r_1\}$ is a supporting hyperplane at $b$, so $S\subset B_R \cap \{ x : x_k < r_1\}$. We claim that 
\begin{equation}\label{fs2} 
\{ (y', y_k)\in \R^{k-1}\times \R  : |y'|< \frac 1{2R}, \ 0<y_k < \frac 1{2r_1}\} \subset S^\circ.
\end{equation}
This is clear, since if $y = (y', y_n)$ belongs to the set on the left, then one readily checks that $x\cdot y\le 1$ for all $x\in S \subset B_R \cap \{ x : x_k<r_1\}$, proving \eqref{fs2}. It follows that 
\[
c > |S^\circ| \ge
\frac {|B^{k-1}_1|}{2r_1(2R)^{k-1} } \ge
\frac {|B^{k-1}_1|}{2r(2R)^{k-1}} .
\]
This cannot happen if $r \le \frac {|B^{k-1}_1|}{2c(2R)^{k-1}}$. So for such $r$, it must be the case that $B_r\subset S$, proving \eqref{fs1}.

Now assume that there exists $A>0$ such that $\oO(\delta)\ge A \delta^{1/n}$, and fix a sequence $a_j\in \Omega$ and $\nu_j\in N(a_j)$ 
such that $\dom(a_j) := \delta_j\to 0$, and $|S^\circ(a_j, \nu_j)| < n A^{-n}$. The existence of such sequences follows directly from Theorem \ref{T2}. Upon passing to subsequences (still labelled $a_j, \nu_j, \delta_j$) we may assume that $a_j\to b\in \pp \Omega$. After a translation and a rotation, we may assume that $b=0$ and $\Omega \subset \{ x\in \R^n : x_n>0\}$.

Appealing to \eqref{fs1} with $k=n-1$, we find that $B_r\subset S(a_j, \nu_j)$ with 
%$r = |B^{n-2}_1| A^n /2^{n-1}n R^{n-2}$.
$r = \frac{|B^{n-2}_1| A^n }{2^{n-1}n R^{n-2}}$.
Then the definition of $S(a_j,\nu_j)$ implies that 
\[
\{ a_j + x : x\cdot \nu_j=0, |x|< r\}\subset \Omega \subset \{ x : x_n>0\}
\quad\mbox{for every $j$.}
\]
This implies that $\nu_j\to - e_n $ as $j\to \infty$ and $a_j\to b=0$. Then 
\[
\{ a_j + x : x\cdot \nu_j=0, |x|< r\}  \ \ \longrightarrow  \ \ \{ x : x\cdot- e_n=0, |x|<r\} = B^{n-1}_r\times \{0\}
\]
as $j\to \infty$, in the Hausdorff distance. It follows that
$B^{n-1}_r\times \{0\}\subset \bar \Omega$ and hence, since $\Omega \subset \{x : x_n>0\}$, we conclude that
$B^{n-1}_r\times \{0\}\subset \pp \Omega$. Thus we have found a flat spot.

We omit the proof that if $\pp \Omega$ has a flat spot, then $\oO(\delta) \ge c\delta^{1/n}$ for some $c$, which is a very direct consequence of Theorem \ref{T2}.
\end{proof}

Finally, we present the proof of a fact already stated in the introduction.

\medskip
\begin{lemma}\label{lem:ellipse.section}
If $x\in \pp\Omega$ is a point at which $\pp \Omega$ is twice
differentiable, with Gaussian curvature $\kappa$, and if $\nu$ is the outer unit normal at $x$, then
\begin{equation}\label{section.poscurv}
|S^\circ(x-\delta\nu, \nu)| = \frac{\sqrt \kappa |B^{n-1}_1|}{(2\delta)^{(n-1)/2}} (1+o(1)) \qquad\mbox{ as }\delta\to 0.
\end{equation}
\end{lemma}

\begin{proof} 
Choosing coordinates so that $x=0$ and $\nu = - e_n$, we find that locally near $0$, $\Omega$ has the form 
$\{ x = (x', x_n)\in \R^{n-1}\times \R: x_n > h(x')\}$ for $h$ such that $h(x') = \frac 12 x'\cdot Q x'(1+o(1))$ as $x'\to 0$, with $\det Q = \kappa$. From there the definitions imply that 
\[
S(x-\delta\nu, \nu) = \{ x'\in \R^{n-1} : h(x') < \delta \} \times \{0\} . 
\]
The expansion of $h$ for small $x'$ implies that for any $\ep>0$, there exists $\delta_0>0$ such that
if $0<\delta<\delta_0$, then
\begin{align*}
\{x'\in \R^{n-1} : x'\cdot Q x' < 2\delta(1-\ep) \} 
&\subset 
\{ x'\in \R^{n-1} : h(x') < \delta \} \\
&\subset 
\{x'\in \R^{n-1} : x'\cdot Q x' < 2\delta(1+\ep) \}.
\end{align*}
Since the ellipse $\{ x' : x'\cdot Q x' < r^2\}$ has volume $r^{n-1}|B^{n-1}_1|/\sqrt{\det Q}$, we deduce \eqref{section.poscurv} 
from the standard fact that  $|E| \, |E^\circ| = |B^{n-1}_1|^2$ for any ellipse $E$ in $\R^{n-1}$, a consequence of affine invariance.
\end{proof}

\section{Proof of Theorem \ref{T1}}

The proof of Theorem \ref{T1} is distributed among Propositions \ref{lem:liminf.oO}, \ref{lem:n=2} and \ref{lem:nge3}. Before starting their proofs, we give a preliminary lemma.

\begin{lemma}\label{lem:ellipse}
Let $E\subset \R^n$ denote the ellipsoid
\[
E := \{ x\in \R^n : \frac {x_1^2}{\ell_1^2} +\cdots + \frac {x_n^2}{\ell_n^2} \le 1 \}
\]
and let $a := (0,\ldots, 0, -\alpha)$ for some $\alpha\in [0,\ell_n)$, and let $p = (0,\ldots,0, -\ell_n)\in \pp E$. If $\alpha$ is close enough to $\ell_n$, then
\[
f_E(a) = \frac{ \sqrt{\kappa(p)} |B^n_1|}{ [d_{E}(a) (2 - \ell_n^{-1} d_{ E}(a))]^{(n+1)/2} } .
\]
\end{lemma}

\begin{proof}
We recall that if $S\subset \R^n$ is a convex set, the support function $\sigma_S$ is defined by
\[
\sigma_S(y) = \sup_{x\in S}x\cdot y.
\]
It is rather clear from the definitions that
\[
S^\circ = \{ y\in \R^n : \sigma_S(y)\le 1\}\qquad\quad
\sigma_{S-a}(y) = \sigma_S(y) - a\cdot y,\qquad\quad
\sigma_{B^n_1}(y) = |y|.
\]

{\em Step 1}. 
Let $B$ denote the unit ball $B^n_1$. Then using the properties of the support function noted above,
\[
f_B(a) = |(B-a)^\circ| = |\{ y\in \R^n : \sigma_{B-a}(y) \le 1\}| = |\{ y\in \R^n  :  |y|-a\cdot y \le 1\}|.
\]
For $a$ as above, by writing $y = (y', y_n)\in \R^{n-1}\times \R$,  by squaring both sides, completing a square and rearranging, we find that
\[
|y| \le 1 + a\cdot y = 1-\alpha y_n
\quad \iff \quad (1-\alpha^2) |y'|^2 + (1-\alpha^2)^2(y_n + \frac{\alpha}{1-\alpha^2})^2   \le 1  .
\]
The inequality on the right defines an ellipsoid whose volume is easily found, yielding
\[
f_B(a) = |(B-a)^\circ| = (1-\alpha^2)^{-(n+1)/2}|B^n_1|.
\]
Since $d_{\partial B}(a) = 1-\alpha$, we can rewrite this as
\[
f_B(a) = \frac{|B^n_1|}{ [d_{ B}(a)(2 - d_{ B}(a)]^{\frac{n+1}2}}
\]

{\em Step 2}. Now let $E$ denote a general ellipsoid as in the statement of the theorem. 
Noting that
\[
E = MB \quad\mbox{ for }M = \mbox{diag}(\ell_1,\ldots, \ell_n).
\]
we find from Lemma \ref{lem:affine} that
\[
f_E(a) = f_{MB}(a)  = (\ell_1 \cdots \ell_n)^{-1} f_B(M^{-1}a),
\]
Also, since $M^{-1}a = (0,\ldots, 0, -\alpha/\ell_n)$, we see that $d_B(M^{-1}a) = 1-\frac \alpha{\ell_n} = \frac{\ell_n-\alpha}{\ell_n}$.. Substituting into the above formula, we obtain
\[
f_E(a) 
%= \frac 1{\ell_1\cdots\ell_n} \frac{|B^n_1|}{ [ \ell_n^{-1}d_{ E}(a)(2 - \ell_n^{-1}d_{ E}(a)]^{\frac{n+1}2}}
= \frac {\ell_n^{\frac{n-1}2}}{\ell_1\cdots\ell_{n-1}} \frac{|B^n_1|}{ [ (
\ell_n-\alpha) (2 - \frac{\ell_n-\alpha}{\ell_n})]^{\frac{n+1}2}}.
\]
It is clear that if $\alpha$ is close enough to $\ell_n$, then $\ell_n-\alpha = d_{ E}(a)$.
So to complete the proof we must show that $\kappa(p) = \ell_n^{n-1}/(\ell_1^2\cdots\ell_{n-1}^2)$.
This is an easy computation. Near $p$, we write $\pp E$ as the graph
\[
x_n = g(x'), \qquad g(x') = -\ell_n\left(1 - ( \frac {x_1^2}{\ell_1^2} +\cdots + \frac {x_n^2}{\ell_{n-1}^2})\right)^{1/2}.
\]
We compute that $Dg(0)=0$ and  $D^2 g(0) = \ell_n \mbox{diag}(\ell_1^{-2}, \ldots , \ell_{n-1}^{-2})$, and it follows that
$\kappa(p) = \det D^2 g(0) = \ell_n^{n-1}/(\ell_1^{2} \cdots  \ell_{n-1}^{2})$ as claimed.
\end{proof}

\begin{proposition}
    \label{lem:liminf.oO}
Let $\Omega\subset \R^n$ be a bounded, convex open subset of $\R^n$ for $n\ge 2$, and assume that \eqref{kmin} holds.
Then
\begin{equation}\label{liminf.oO}
\liminf_{\delta\to 0}  \frac{\oO(\delta)}{\delta^{\frac {n+1}{2n}}} \ge 
 \left(\frac{2^{(n+1)/2}}{|B^n_1| \sqrt{\kappa_0}}
\right)^{1/n}.
\end{equation}
\end{proposition}

\begin{proof}
Given $\ep>0$, choose a point $p_\ep\in \pp\Omega$ at which $\pp \Omega$ is twice differentiable and
$\kappa(p_\ep) < (1+\ep)\kappa_0$. We may assume after a translation and a rotation that $p_\ep=0$ and that there exists $r>0$ such that in a neighborhood of $p_\ep$, 
\begin{equation}\label{graph}
\Omega\cap B^n_r  = \{ (x', x_n)\in B^n_r : x_n > g(x') \}
\end{equation}
for a convex $g$ such that
\begin{equation}\label{g.expand}
g(x') = \frac 12 \sum_{j=1}^{n-1}\lambda_j x_j^2  + o(|x'|^2) \quad\mbox{ as } |x'|\to 0, 
\end{equation}
with
\[
\prod_{j=1}^{n-1} \lambda_j = \kappa(p_\ep) = \kappa(0), \qquad \lambda_j>0\mbox{ for all }j.
\] 
Now let $E_\ep$ be the ellipse
\[
E_\ep = \{ x\in \R^n : \frac {x_1^2}{\ell_1^2} + \cdots+ \frac {x_{n-1}^2}{\ell_{n-1}^2} + \frac{(x_n - \ell_n)^2}{\ell_n^2} < 1\}
\]
for 
\[
\ell_n = \eta\mbox{ to be chosen}, \quad \ell_j = \left( \frac{\ell_n}{\lambda_j(1+\ep)}\right)^{1/2}.
\]
We claim that
\begin{equation}\label{small.ellipse}
\mbox{if $\eta>0$ is taken to be small enough, then } E_\ep \subset \Omega.
\end{equation}
In view of \eqref{graph}, it suffices to show that if $\eta$ is small enough and $x = (x',x_n)\in E_\ep$, then
$x\in B_r$ and $x_n > g(x')$.

It is clear that there exists $\eta_0>0$ such that  $E_\ep \subset B^n_r$ whenever $0<\eta<\eta_0$.

Second, we use the concavity of the square root to see that
\begin{align*}
x\in E_\ep 
&\quad\Longrightarrow \quad 
\frac{(\ell_n-x_n)^2}{\ell_n^2}\le 1 - \sum_{j=1}^{n-1} \frac {x_j^2}{\ell_j^2} \\
&\quad\Longrightarrow \quad 
\frac{(\ell_n-x_n)}{\ell_n}\le  \left(1 - \sum_{j=1}^{n-1} \frac {x_j^2}{\ell_j^2}  \right)^{1/2} \le 1 - \frac 12\sum_{j=1}^{n-1} \frac {x_j^2}{\ell_j^2}  \\
&\quad\Longrightarrow \quad 
x_n \ge \frac 12 \sum_{j=1}^{n-1} \frac{\ell_n}{\ell_j^2} x_j^2 = \frac {(1+\ep)}2 \sum_{j=1}^{n-1}\lambda_j x_j^2.
 \end{align*}
Note also that the definitions imply that $|x'| < C \sqrt\eta$ for $(x',x_n)\in E_\ep$, so the above inequality and  \eqref{g.expand} implies that there exists $\eta_1\in (0,\eta_0)$ such that if $0<\eta<\eta_1$, then $x_n > g(x')$, completing the proof of \eqref{small.ellipse}.

We henceforth fix $\eta<\eta_1$.
Let $a_\delta = (0,\ldots, 0,\delta)$ for $\delta< \ell_n$, and note that $a_\delta\in E_\ep$ when 
$\delta<2\ell_n$. Note also that $\dom(a_\delta) = d_{ E_\ep}(a_\delta)= \delta$ for all
sufficiently small $\delta>0$.

It follows from Lemma \ref{lem:compare} that $\fom(a_\delta)\le f_{E_\ep}(a_\delta)
(a)$, so we use (an easy modification of) Lemma \ref{lem:ellipse} to conclude
\[
\lim_{\delta\searrow 0} \delta^{\frac{n+1}2}\fom(a_\delta) \le
\lim_{\delta\searrow  0} \delta^{\frac{n+1}2}f_{E_\ep}(a_\delta) = \lim_{\delta\searrow 0} \frac {\sqrt{\kappa_{E_\ep}(0)}|B^n_1|}{(2- \frac \delta {\ell_n})^{(n+1)/2}} =
\frac {\sqrt{\kappa_{E_\ep}(0)}|B^n_1|}{2^{(n+1)/2}},
\]
where $\kappa_{E_\ep}(0)$ denotes the curvature of $\pp E_\ep$ at $0$, which is
\[
\kappa_{E_\ep}(0) = (1+\ep)^{n-1}\kappa(0) < (1+\ep)^n \kappa_0.
\]
Applying Proposition \ref{prop:fom}, we find that
\[
\liminf_{\delta\searrow 0} \delta^{-\frac{n+1}{2n}}  \oO(\delta) >  \frac 1{\sqrt{1+\ep}} \left( 
 \frac{2^{(n+1)/2}}{ {\sqrt \kappa_0}|B^n_1|}
 \right)^{1/n}.
\]
Since $\ep>0$ was arbitrary,  conclusion \eqref{liminf.oO} follows.
\end{proof}

\begin{proposition}  \label{lem:n=2}
Let $\Omega$ be a bounded, open, convex subset of $\R^2$ satisfying \eqref{kmin},
then
\[
\sup_{\delta>0} \frac{\oO(\delta)}{\delta^{ 3/4} }= \left( \frac
{2^{3/2}}{\sqrt{\kappa_0} \pi}
\right)^{1/2}.
\]
\end{proposition}

\begin{proof}
We will show that for any $a\in \Omega$,
\begin{equation}\label{inf.fom}
\fom(a) \ge \frac 1{\dom(a)^{ \frac 32} } \frac{ \sqrt{\kappa_0} \pi}{2^{3/2} }
\end{equation}
In view of Proposition \ref{prop:fom}, this implies that
\[
\sup_{\delta>0} \frac{\oO(\delta)}{\delta^{ \frac 34 } } \le \left( \frac {2^{3/2}}
{\sqrt{\kappa_0} \pi}\right)^{1/2}.
\]
This will complete the proof of the Proposition, as the opposite inequality follows from Proposition \ref{lem:liminf.oO}.

To prove \eqref{inf.fom}, fix any $a\in \Omega$, and let $b\in \pp \Omega$ be a point such that
$\dom(a) = |a-b| =: \delta$, not necessarily small. After a rotation and a translation, we may assume that $b = (0,0)$ and $a = (0,\delta)$. Clearly $B_\delta(a)\subset\Omega$, and $0\in \pp B_\delta(a)\cap \pp\Omega$. From these facts and the convexity of $\Omega$ one easily sees that 
$\Omega\subset \{ (x_1,x_2) : x_2>0\}$
and
\[
\pp \Omega\cap [(-\delta,\delta)\times (0,\delta)] \ \ 
\subset \ \  \left\{( x_1,x_2) : 0< x_2 \le \delta - \sqrt{ \delta^2 - x_1^2}\right\}. 
\]
Let $I = \{ x_1 : (x_1,x_2)\in \Omega\mbox{ for some }x_2\}$ be the projection of $\Omega$ onto the $x_1$-axis. Then writing the lower part of $\pp\Omega$ as the graph of a function $g:I\to \R$, we have
\[
g\mbox{ is convex}, \qquad 0\le g(x_1) \le \delta - \sqrt{\delta^2-x_1^2},\qquad \Omega\subset \{(x_1,x_2) : x_1\in I, x_2> g(x_1)\}.
\]
Note that $g$ is differentiable at $x_1=0$, with $g'(0)=0$.

We now claim that 
\begin{equation}\label{n=2.c1}
g(x_1) \ge \frac {\kappa_0}2 x_1^2 \quad\mbox{ for }x_1\in I, \quad \mbox{ and thus }\Omega \subset D :=  \{ (x_1,x_2) : x_2>  \kappa_0  x_1^2/2 \}.
\end{equation}
We will prove this for $x_1>0$; the argument for $x_1<0$ is basically identical.
Let us write $S := \{ x_1\in I : g\mbox{ is twice differentiable at }x_1\}$. To prove \eqref{n=2.c1}, we recall assumption \eqref{kmin}, which implies that
\[
\frac{g''(x_1)}{(1+(g'(x_1)^2)^{3/2}} \ge \kappa_0\qquad\mbox{ for every }x_1\in S
\]
This clearly implies that $g''(x_1)\ge \kappa_0$ in $S$. Since $g$ is convex, $g'$ is increasing function, so
for positive $x_1\in I$, elementary real analysis yields
\[
g'(x_1) = g'(x_1) - g'(0) \ge \int_{\{t\in S : 0<t< x_1\}} g''(t) dt  \ge \kappa_0 x_1.
\]
Since $g'$ is locally Lipschitz, we obtain \eqref{n=2.c1} by integrating again.

In view of \eqref{n=2.c1} and  Lemma \ref{lem:compare}, in order to prove \eqref{inf.fom}, it suffices to show that
\begin{equation}\label{n=2.c2}
f_D(a) = \frac {\sqrt{\kappa_0}\pi}{(2\delta)^{3/2}}.
\end{equation}
This is a straightforward computation. First, let 
\[
M := \left(\begin{array}{cc}  \sqrt{\frac {\kappa_0}{2\delta} }&0\\
0&\frac 1 \delta\end{array}\right)
\]
Then $\widetilde D := MD = \{ (x_1, x_2) : x_2> x_1^2 \}$ and $Ma = (0,1) = e_2$. By Lemma
\ref{lem:affine}, 
\[
f_D(a) = \sqrt{\frac{ \kappa_0}{2\delta^{3/2} } }\ f_{\widetilde D}(e_2)
\]
And
\[
f_{\widetilde D}(e_2) = |(\widetilde D - e_2)^\circ| = |\{ y\in \R^2 : \sigma_{\widetilde D - e_2}(y) \le 1\}|.
\]
Given $y\in \R^2$, we compute $\sigma_{\widetilde D - e_2}(y)$ by attempting to find $x$ that maximizes $x\mapsto y\cdot (x-e_2)$ subject to the constraint that $x\in \widetilde D$.
It is clear that a maximum can only occur for $x\in \pp \widetilde D$, so one can use Lagrange multipliers to find that
\[ 
\sigma_{\widetilde D - e_2}(y_1,y_2) = \begin{cases}
+\infty&\mbox{ if }y_2\ge 0\\
-\frac {y_1^2}{4y_2} - y_2 &\mbox{ if }y_2<0. 
\end{cases}
\]
It easily follows that
\[
(\widetilde D-e_2)^\circ = \left\{ (y_1, y_2) : y_1^2 + 4\left(y_2-\frac 12\right)^2 \le 1\right\}, \qquad\mbox{ and hence }\quad f_{\widetilde D}(e_2) = \frac \pi 2,
\]
concluding the proof of \eqref{n=2.c2}.
\end{proof}

The remaining assertion of Theorem \ref{T1} is contained in the following proposition. The idea of the proof is to
approximate $\pp \Omega$ from the outside, locally, by a quadratic. We need to be able to do this in a uniform way, and having done so, to extract information about $\pp u_a(\Omega)$ from its quadratic approximation when this approximation is only local.

\begin{proposition}\label{lem:nge3}
If $n\ge 3$ and $\pp\Omega$ is $C^2$, then 
\begin{equation}\label{eq:nge3}
\lim_{\delta \searrow 0} \frac{\oO(\delta)}{\delta^{(n+1)/2n} } =  
 \left(\frac{2^{(n+1)/2}}{|B^n_1| \sqrt{\kappa_0}}
\right)^{1/n}
\end{equation}
\end{proposition}

\begin{proof}
In view of Proposition \ref{lem:liminf.oO} and Proposition \ref{prop:fom}, we only need to prove that
\begin{equation}\label{fom.lbd.kpos}
\liminf_{\delta\searrow 0} \inf_{\dom(a)=\delta} \delta^{(n+1)/2} \fom(a)  \ge 
\frac{|B^n_1| \sqrt{\kappa_0}}
{2^{(n+1)/2}}.
\end{equation}

{\em Step 1}. We first claim that for any $\ep_1>0$, there exists $r_0>0$ such that for any $b  \in \pp\Omega$, 
there exists a rigid motion $S$ (that is, the composition of a rotation and a translation)
and a convex $C^2$ function $g:B^{n-1}_{r_0}\to [0,\infty)$ such that
\[
S(0) = b, \qquad S(\{ (x', g(x')): |x'|<r_0 \}) \subset \pp \Omega,
\]
and $g(0)=0$, with
\begin{equation}\label{unif.C2}
\|D^2 g(x') - D^2 g(0)\| < \ep_1
\qquad\mbox{ for }x'\in B_{r_0}(0)
\end{equation}
where $\|  \cdot  \|$ denotes the operator norm. 
Informally, this states that $\pp \Omega$ is {\em uniformly} $C^2$. Since $\pp \Omega$ is $C^2$ and compact, on some level this is clear, but we provide some details nonetheless. Our proof of \eqref{unif.C2} will also show that there exist positive $\Lambda_{min}\le \Lambda_{max}$, independent of $b\in \pp\Omega$, such that
\begin{equation}\label{LllL}
\Lambda_{min} \le \lambda_1\le \cdots \le \lambda_{n-1}\le\Lambda_{max}, \qquad
\{ \lambda_j \} = \mbox{ eigenvalues of }D^2g(0).
\end{equation}

First, the compactness of $\pp \Omega$ implies that  there exists $R>0$, 
$J\ge 2$
and
\begin{itemize}
\item maps $S_j: \R^n\to \R^n$ for $j=1,\ldots, J$,  each one a rigid motion (the composition of a translation and a rotation),
\item  $C^2$ functions $h_j:B_{4R}^{n-1}\to [0,\infty)$  for $j=1,\ldots, J$,
\end{itemize}
such that $|\nabla h_j| \le \frac 14$ in 
$B_{4R}^{n-1}$
and 
\[
\pp \Omega = \cup_{j=1}^J  U_j, \qquad U_j :=S_j \left( \{ (x', h_j(x')) : x'\in B^{n-1}_R \} \right).
\]
For every $j$,  clearly $h_j, D h_j$ and $D^2 h_j$ are uniformly continuous on $B_{3R}^{n-1}$, and there are only finitely many of these functions, so there
exists a common {\em $C^2$ modulus of continuity} for $\{ h_j\}_{j=1}^J$ on $B_{3R}^{n-1}$, by which we mean a 
continuous, increasing function $\mu:[0,\infty)\to[0,\infty)$ such that $\mu(0)=0$ and 
\begin{equation}\label{C2mod}
|h_j(x)-h_j(y)| + |Dh_j(x)-Dh_j(y)| + |D^2h_j(x)-D^2 h_j(y)| \le \mu(|x-y|)
\end{equation}
for all $x$ and $y$ in $B_{3R}^{n-1}$ and $j=1,\ldots, J$.

Now consider any $b\in \pp \Omega$. Fix some $\bar x'\in B^{n-1}_R(0)$ and $j\in  \{1,\ldots, J\}$ such that $b = S_j(\bar x', h_j(\bar x'))$. By a further translation we may send $(\bar x', h_j(\bar x'))$ to the origin in $\R^n$, then by a rotation in $\R^n$ we can arrange that part of the (translated and rotated) graph of $h_j$ over $B_{2R}(\bar x')$ can be written as the graph over some ball $B_{r}^{n-1}$ of a convex $C^2$ function $g:B^{n-1}_r\to [0,\infty)$ such that $g(0)=0$ and $\nabla g(0)=0$. By applying  Lemma \ref{lem:rotate} below to $h(x') = h_j(x' - \bar x') - h_j(\bar x')$, we find that independent of the choice of $b$, the domain of the resulting function $g$ can be taken to be $B^{n-1}_R$, and $g$ has a $C^2$ modulus of continuity in $B_R^{n-1}$ that can be estimated solely in terms of $\mu$ from \eqref{C2mod}. This proves \eqref{unif.C2}. Similarly, \eqref{LllL} follows from \eqref{d2gy0} in Lemma \ref{lem:rotate} which we prove below,  together with the fact that for every $j\in \{1,\ldots,J\}$ the eigenvalues of 
$D^2 h_j$ are bounded away from $0$ in $B^{n-1}_R$, being positive on $B^{n-1}_{2R}$.

\medskip

{\em Step 2}. We now prove \eqref{fom.lbd.kpos}. 

{\em Step 2.1:  Normalization and approximation by a quadratic}.
Because $\pp \Omega$ is $C^2$ and compact, there exists $\delta_0>0$ such that if $\delta:= \dom(a)  <\delta_0$, then there is a unique  $b\in \pp\Omega$ such that $\dom(a)=|a-b|$. Fix some such $a$ and $b$. 
In view of Step 1, we may assume after a rigid motion that $a = \delta e_n = (0,\ldots, 0,\delta)$ and 
$b=0$, and that there is a  nonnegative convex function $g$, vanishing at $x'=0$, such that 
\eqref{unif.C2} holds for some $\ep_1$ and $r_0(\ep_1)$ to be specified in a moment,  and
with $\{(x', g(x')) : |x'|<r_0\}$ contained in the (rotated and translated) $\pp \Omega$.

Fix $\ep>0$.  Using \eqref{unif.C2} for a suitably small choice of $\ep_1$ and calculus, we can guarantee that
\begin{equation}\label{choose.ep1}
\left. \begin{aligned}
g(x') > \frac {(1-\ep)}2 x'\cdot Qx' &\\ 
|Dg(x') - Q x'| < \ep |x'| &
\end{aligned}\right\}
\  \mbox{ in }B_{r_0}^{n-1}, \quad\mbox{ for }Q = D^2g(0).
\end{equation}

Let $M:= Q^{1/2}$,  the positive definite symmetric square root of $Q$, and define
\[
\widetilde \Omega := \left\{ \left(\frac {M x'}{\sqrt \delta}, \frac{x_n}\delta \right) : (x', x_n)\in \Omega\right\}, \qquad
\tilde g(y') := \frac 1 \delta g(\sqrt \delta M^{-1} y').
\]
Due to \eqref{LllL}, the eigenvalues of $M$ are bounded between
$ \Lambda_{min}^{1/2}$ and $\Lambda_{max}^{1/2}$.
Thus
\begin{equation}\label{rdelta.def}
|y'|< r_\delta :=  \frac{r_0}{\sqrt{\delta\Lambda_{min}}} \qquad\Longrightarrow \quad |\sqrt \delta M^{-1} y'| <  r_0 .
\end{equation}
The definition of $\widetilde g$ is chosen so that
\[
\{ (y', \tilde g(y')) : |y'| < r_\delta \} \ \subset \  \pp \widetilde \Omega.
\]
Hypothesis \eqref{kmin} implies that $\det Q = (\mbox{curvature of }\pp\Omega\mbox{ at }b=0) \ge \kappa_0$,
so we deduce from Lemma \ref{lem:affine} that
\begin{equation}\label{ftildeO}
\delta^{(n+1)/2}\fom(a) = |\! \det M| \, f_{\widetilde \Omega}(\tilde a) \ge   \sqrt \kappa_0\ f_{\widetilde \Omega}(\tilde a) 
\qquad\mbox{ for } \ \ \tilde a = \frac a \delta = e_n.
\end{equation}
In addition, using \eqref{rdelta.def} and requiring that
$\ep< \frac 12 \Lambda_{min}$, we can translate properties \eqref{choose.ep1} into statements about $\tilde g$. It follows that$(y',\tilde g(y'))\in \pp \widetilde \Omega$ and
\begin{equation}\label{tilde.g}
\left. \begin{aligned}
\tilde g(y') > \frac {(1-\ep)}2 |y'|^2 & \\ 
|D\tilde g(y') - y'| < \ep \Lambda_{min}^{-1} |y'| <\frac 12 |y'| &
\end{aligned}\right\} \mbox{ in }B_{r_\delta}^{n-1}.
\end{equation}
Fix $z'$ such that $|z'|< r_\delta/2$ and let $\Phi(y') := z'+y' - D\tilde g(y')$. 
Then the second inequality in \eqref{tilde.g} implies that $\Phi$ maps
$\bar B_{r_2}^{n-1}$ to itself, for $r_2 =  2|z'| < r_\delta$.
Thus the Brouwer Fixed Point Theorem implies that $\Phi$ has
a fixed point $y'$ in $B_{r_2}^{n-1}$. But $\Phi(y') =y'$ exactly when $D\tilde g(y') = z'$.
It follows that 
\begin{equation}\label{dtildeg.image}
B^{n-1}_{r_\delta/2} \subset \{ D\tilde g(y') :  |y'|< |B^{n-1}_{r_\delta}| \}.
\end{equation}

{\em Step 2.2: finding a large subset of $\pp u_{\widetilde \Omega, a}(\tilde \Omega)$}.

We will write $\tilde u_{\tilde a} := u_{\widetilde \Omega, \tilde a}$. We next will show that
\begin{equation}\label{pptildeu}
E_{\ep,\delta} := \left\{ s( z', -1) :|z'|<  \frac {r_\delta}2, \ \  0\le s \le \left( 1 + \frac 1{2(1-\ep)} |z'|^2 \right)^{-1}  \right \} \subset \pp \tilde u_{\tilde a}(\widetilde \Omega).
\end{equation}
To see this, fix $z'\in B^{n-1}_{r_\delta/2}$, and using \eqref{dtildeg.image}, find $y'\in B^{n-1}_{r_\delta}$ such that $\nabla \tilde g(y') = z'$.
Let 
\[
\ell_{y'}(x) :=  \frac{ (\nabla \tilde g(y'), -1) \cdot(x - (y', \tilde g(y')))}
{1 + y'\cdot \nabla \tilde g(y')-  \tilde g(y')}
\]
We claim that $\ell_{y'}$ is a supporting hyperplane to the graph of $\tilde u_{\tilde a}$ at $\tilde a =e_n$.
We must show that $\ell_{y'}(\tilde a) = \tilde u_{\tilde a}(\tilde a) = -1$, which follows directly from the definition, and  that $\ell_a \le \tilde u_{\tilde a}$ in $\widetilde \Omega$.
Since both $\ell_a$ and $\tilde u_{\tilde a}$ are linear on line segments connecting $\pp \widetilde \Omega$ to $\tilde a$, it suffices to check that $\ell_a \le \tilde u_{\tilde a} = 0$ on $\pp \widetilde \Omega$.
This follows from noting that $\ell_{y'}$ vanishes exactly on the hyperplane 
$\{ x\in \R^n : \nu(y) \cdot (x -  y) = 0\}$, where $y = (y', \tilde g(y'))\in \pp\widetilde\Omega$ and $\nu(y)$ is the outer unit normal to $\pp \widetilde\Omega$ at $y$. This is a supporting hyperplane to $\pp \widetilde \Omega$, so $\ell_{y'}$ does not change sign in $\widetilde \Omega$. Since $\ell_{y'}(\tilde a)<0$, the claim follows.

Thus $\nabla \ell_{y'}(\tilde a)\in \pp \tilde u_{\tilde a}(\tilde a)$. Since it is clear that $0\in \pp\tilde u_{\tilde a}(\tilde a)$ and $\pp\tilde u_{\tilde a}(\tilde a)$ is convex, it follows that the segment  $\{ s\nabla \ell_{y'}(\tilde a) : 0\le s \le 1\} \subset \pp \tilde u_{\tilde a}(\tilde a)$, that is, 
\begin{equation}\label{combining}
\left\{ s(\nabla\tilde g(y') , -1) \ : \  0 \le s \le \frac 1{1+y'\cdot \nabla\tilde g(y')-\tilde g(y')} \right \} \subset \pp \tilde u_{\tilde a}(a).
\end{equation}
However, recalling that $\nabla\tilde g(y') = z'$, and using \eqref{tilde.g} and elementary inequalities,
\[
y'\cdot \nabla \tilde g(y') - \tilde g(y') \le y'\cdot z' - \frac {(1-\ep)}2|y'|^2 \le \frac 1{2(1-\ep)}|z'|^2
\]
We deduce \eqref{pptildeu} by combining this and \eqref{combining}.

{\em Step 2.3: conclusion of proof}.
Since $a$ was an arbitrary point such that $\dom(a)=\delta$, it follows from \eqref{pptildeu} and \eqref{ftildeO} that for all sufficiently small $\delta$,
\begin{equation}\label{previous}
\inf_{\dom(a)=\delta}  \delta^{(n+1)/2}\fom(a) 
\ge \sqrt{\kappa_0} |E_{\ep,\delta}|.
\end{equation}
We will show that
\begin{equation}\label{volA.est}
\lim_{\delta\to 0}|E_{\ep,\delta}| = (1-\ep)^{(n-1)/2}\frac {|B^n_1|}{2^{ (n+1)/2}}
\end{equation}
Since $\ep>0$ is arbitrary, this and \eqref{previous} imply
\eqref{fom.lbd.kpos}, and thus complete the proof of the Proposition.

To establish \eqref{volA.est}, note first that $E_{\ep,\delta}$ forms an increasing family of sets as $\delta\searrow 0$. Thus the Monotone Convergence Theorem implies that
$\lim_{\delta \to 0}|E_{\ep,\delta}| = |E_{\ep,0}|$, 
for
\[
E_{\ep,0} := \cup_{\delta>0} E_{\ep,\delta} = \left\{ s(p', -1) \,:\ p'\in \R^{n-1},  \  0 \le s \le
\left(1 + \frac 1{2(1-\ep)} |p'|^2 \right)^{-1}  \right\}.
\]
We claim that in fact
\begin{equation}\label{EcalE}
E_{\ep,0} = \left\{ (q', q_n) \in \R^n : \frac 2{(1-\ep)}|q'|^2 + 4(q_n  +\frac 12)^2 \le 1 \right\} =: \mathcal E_{\ep,0}.
\end{equation}
Indeed, both $\mathcal E_{\ep,0}$ and $E_{\ep,0}$ are contained in the set $\{ 0\} \cup \{(q', q_n) : q_n < 0\}$. It is clear that the origin belongs to both sets. Any point with $(q', q_n)$ with $q_n<0$ can be written
\begin{equation}\label{repq}
(q', q_n) =  \frac{ t(p',-1)}{1 +  \frac 1{2(1-\ep)} |p'|^2 } \qquad\mbox{ for some $p'\in \R^{n-1}$ and }t>0.
\end{equation}
Then
\[
|q'|^2 = t^2 \frac{|p'|^2}{(1 +  \frac 1{2(1-\ep)} |p'|^2)^2}
\]
and
\[
4(q_n+\frac 12)^2 = (2q_n+1)^2 = t^2  \frac{\left( -1+ \frac 1{2(1-\ep)} |p'|^2\right)^2}{(1 +  \frac 1{2(1-\ep)} |p'|^2)^2}
\]
from which we see  that 
\[
\frac 2{(1-\ep)}|q'|^2 + 4(q_n  +\frac 12)^2 = t^2 .
\]
This implies \eqref{EcalE}, since for $q_n<0$, 
\[
(q', q_n)\in E_{\ep,0} \quad\Longleftrightarrow\quad 0 <t\le 1\mbox{ in }\eqref{repq} 
\quad\Longleftrightarrow\quad (q',q_n)\in \mathcal E_{\ep,0}.
\]
And it is clear that $|\mathcal E_{\ep,0}|= \frac{(1-\ep)^{(n-1)/2}|B^n_1|}{2^{(n+1)/2}}$, proving \eqref{volA.est}.\end{proof}

Finally, we establish the lemma used above.  
In the proof, we find it helpful to write points $x = (x', x_n)\in \R^n = \R^{n-1}\times \R$ as column vectors
$\binom {x'}{x_n}$.

\begin{lemma}\label{lem:rotate}
Let $h\in C^2(B^{n-1}_{2R})$ satisfy $h(0)= 0$ and $\| \nabla h\|_{L^\infty} \le \frac 14$.

Then there exists a rotation $S\in SO(n)$ and a function $g\in C^2(B^{n-1}_R)$ such that 
$\nabla g(0)=0$, 
\begin{equation}\label{rotate1} 
\left\{ \binom{y'}{ g(y')} : |y'| < R \right\} \subset \left\{ S\binom{x'}{ h(x')} : |x'|<20R/11 \right\},
\end{equation}
and such that the $C^2$ modulus of continuity of $g$ can be estimated in terms only of the $C^2$ modulus of continuity of $h$ in $B^{n-1}_{20R/11}$.  Moreover, 
\begin{equation}\label{d2gy0}
g_{y_i y_j}(0) = \begin{cases}
(1+m^2)^{-1/2} h_{x_i x_j}(0) &\mbox{ if }i,j\le n-2\\
(1+m^2)^{-1}h_{x_i x_j}(0) &\mbox{ if }i<j=n-1\mbox{ or }j<i=n-1\\
(1+m^2)^{-3/2} h_{x_i x_j}(0) &\mbox{ if } i=j =n-1\,,\end{cases}
\end{equation}
for some $m \in [-\frac{1}{4},\frac{1}{4}]$, and if $h$ is convex then $g$ is convex.
\end{lemma}

The conclusion of the lemma is a little stronger than we need for the proof or Proposition \ref{lem:nge3}.

\begin{proof}
We may assume by a suitable choice of coordinates that $\nabla h (0) = m e_{n-1}$ for some $m\in [-\frac 14, \frac 14]$, where $e_{n-1}$ is the standard basis vector along the $x_{n-1}$ axis.

We then define (temporarily adopting column vector notation for ease of reading)
\[
S\left(\begin{array}{c}x_1\\ \vdots \\  x_{n-2} \\  x_{n-1} \\ x_n\end{array}\right)
 = \left(\begin{array}{c}x_1\\ \vdots \\ x_{n-2} \\ \alpha x_{n-1}+\beta x_n \\ -\beta x_{n-1} + \alpha x_n\end{array}\right), \qquad \alpha = \frac 1{\sqrt {1+m^2}}, \quad \beta = \frac m{\sqrt {1+m^2}}\ .
\]
Clearly $S\in SO(n)$. Note that $\alpha\ge 4/\sqrt {17} > 4/5$ and $|\beta|\le 1/\sqrt {17} < 1/4$.
We next define $\Phi:B_{2R}^{n-1} \to \R^n$
\[
 S \binom{x'}{h(x')} =: \Phi(x') = \binom {\phi(x')}{\phi^n(x')}
\]
where we will write the components of $\phi$ as $(\phi^1,\ldots, \phi^{n-1})$,
for $\phi$ and $\phi^n$ taking values in $\R^{n-1}$ and $\R$ respectively.
We now define
\[
\psi = \phi^{-1}, \qquad g = \phi^n\circ \psi.
\]
%We will verift below that $\phi$ is injective and that $B^{n-1}_R\subset \phi(B^{n-1}_{3R/2}) $, so that both $\phi$ and $g$ are defined on $B_R^{n-1}$.
We now verify that these functions have the stated properties.

\medskip

{\em Proof that $g$ is well-defined on $B_R^{n-1}$}:
From the definitions we see that 
\begin{equation}\label{Dphi}
\pp_i\phi^j = \delta^j_i \ \ \mbox{ if }j\le n-2,
\qquad
D \phi^{n-1} = \alpha e^{n-1} + \beta D h\ .
%\pp_i \phi^{n-1} = \delta_{i,n-1} \alpha + \beta \pp_i h(x') .
\end{equation}
Thus, writing $\| A\|$ to denote the operator norm of a matrix $A$ and $I_k$ the $k\times k$ identity matrix, we check that
\begin{equation}
\| D\phi - I_{n-1} \| \le |1-\alpha| + |\beta| < \frac 9{20} \qquad\mbox{ everywhere in $B_{2R}^{n-1}$.}
\label{close.to.I}\end{equation}
Now fix any $y'\in B_R^{n-1}$, and define $\Phi(x') = y'+x' -\phi(x')$.
Then for any $x', z'\in B_{2R}^{n-1}$,
\begin{align}
\left|\Phi(x') - \Phi(z')\right|
 =\left| \int_0^1 \frac d{ds}\Phi(s x'+(1-s )z')\,ds\right| 
 &\le  \sup_{B^{n-1}_{2R}} \| D\phi - I_{n-1} \| \, |x'-z'| \nonumber
 \\
 & \le \frac 9{20}|x'-z'|. \label{contraction}
\end{align}
Thus $\Phi$ is a contraction mapping. Note also, that when $y=0$, we find from \eqref{contraction} that $|\phi(x')-x'| \le \frac 9{20}|x'|$. So if $|x'| \le 20R/11$, then 
\[
|\Phi(x')| \le |y'| + |\phi(x')-x'| \le \frac {20R}{11}.
\]
So $\Phi$ maps $B^{n-1}_{20R/11}$ to itself, and the Contraction Mapping Principle thus implies that there is a unique $z'\in B^{n-1}_{20R/11}$ such that $\Phi(z') = z'$, which says exactly that $\phi(z') = y'$.

These facts imply that $\phi^{-1} =\psi$ is well-defined in $B^{n-1}_R$, taking values in $B^{n-1}_{20R/11}$, and hence that $g$ is well-defined in $B^{n-1}_R$ as well.

{\em Proof that \eqref{rotate1} holds}: The definitions imply that if $y' = \phi(x')$, then $g(y') = \phi^n(x')$, and thus
\[
\binom {y'}{g(y')} = \binom{\phi(x')}{\phi^n(x')} = \Phi(x') = S\binom {x'}{h(x')}.
\]
We deduce \eqref{rotate1} from this and remarks above about the range of $\psi$.

{\em Proof that $\nabla g(0)=0$}:
We compute
\begin{equation}\label{dyg}
g_{y_i} =   (\phi^n_{x_k}\circ \psi )\, \psi^k_{y_i}. 
\end{equation}
Since $\phi(0)=0$, it is clear that $\psi(0)=0$. It thus suffices to check that $\nabla \phi^n (0)=0$.
This follows from the choice of $\alpha$ and $\beta$, which guarantees that
\[
\phi^n_{x_k}(0) = \begin{cases} \alpha h_{x_k}(0)  = 0&\mbox{ if }k\le n-2\\
-\beta + \alpha h_{x_{n-1}}(0) = 0&\mbox{ if }k=n-1
\end{cases}
\]

{\em $C^2$ modulus of continuity of $g$:}
By differentiating \eqref{dyg} we obtain
\begin{equation}\label{d2gy}
g_{y_i y_j} = (\phi^n_{x_k x_l}\circ \psi)\psi^k_{y_i}\psi^l_{y_j} + (\phi^n_{x_k}\circ \psi )\psi^k_{y_i y_j}.
\end{equation}
Moreover, 
\[
\psi^i_{y_j} = [(D\phi)^{-1}\circ \psi]^i_j,
\qquad
\psi^i_{y_j y_k} = -  (\phi^m_{x_a x_b}\circ \psi)\, \psi^i_{y_m}\ \psi^a_{y_j}\ \psi^b_{y_k}.
\]
For every $x'\in B_{2R}^{n-1}$, it follows from \eqref{close.to.I} that  $D\phi(x')$ belongs to
the compact set $\{ A\in M^{n-1} : \| A - I\| \le 9/20 \}$ on which the map $A\mapsto A^{-1}$ is smooth and hence bounded and Lipschitz. Hence $D\psi$ is bounded in $B_R^{n-1}$, and $\psi$ is Lipschitz continuous. 
Then elementary estimates show that first $D\psi$ and then $D^2\psi$ have moduli of continuity estimated only in terms of the $C^2$ modulus of continuity of $\phi$, which in turn is controlled by the $C^2$ modulus of continuity of $h$. Then similar arguments show that the $C^2$ modulus of continuity of $g$ can be estimated only in terms of that of $h$. 

{\em Formula for $D^2 g(0)$}. Computing as in our verification that $\nabla g(0)=\nabla \phi^n(0)=0$, and recalling that $\psi(0)=0$, one easily checks that $D^2 \phi^n(0) = \alpha D^2 h(0)$. Similarly, nothing from the definitions that 
$\pp_{n-1}\phi^{n-1}(0) = \alpha+\beta m = 1/\alpha$, one checks that
\[
D\phi(0) = \mbox{diag}(1,\ldots,1,\frac 1 \alpha),
\quad\mbox{  and so  }\quad
D\psi(0) = \mbox{diag}(1,\ldots,1,  \alpha).
\]
We deduce \eqref{d2gy0}  from these facts and \eqref{d2gy}.

{\em Finally}, it is clear that if $h$ is convex then $g$ is convex, as then the graph of $g$ is part of the lower boundary of a convex set.
\end{proof}

\medskip\noindent
{\it Acknowledgements.} This work was partially supported by the Natural Sciences and Engineering Research Council of Canada under
Operating Grant 261955. In addition, the research of Charles Griffin was partially funded by an NSERC Undergraduate Summer Research Award. The authors thank Nam Q. Le for some helpful remarks.

\bibliography{Alexandrov.bib}
\bibliographystyle{acm}

\end{document}